\documentclass[12pt]{article}
\usepackage{amssymb}
\usepackage{mathrsfs}
\makeatletter
 
 \@addtoreset{equation}{section}
\makeatother

\newcommand{\prf}{\noindent{\bf Proof. }}

\newcommand{\rem}{\noindent{\bf Remark. }}

\newcommand{\ff}[1]{%
    {\bf F}_{#1}}

\newcommand{\qed}{\hbox{\rule[-2pt]{5pt}{11pt}}}

\newtheorem{dfn}{Definition}[section]
\newtheorem{thm}[dfn]{Theorem}

\newtheorem{cor}[dfn]{Corollary}

\newtheorem{exam}[dfn]{Example}
\newtheorem{conj}[dfn]{Conjecture}
\newtheorem{prob}[dfn]{Problem}

\begin{document}
%\begin{linenumbers}
%\begin{frontmatter}
\title{Construction of divisible formal weight enumerators and extremal polynomials not satisfying the Riemann hypothesis}
\author{Koji Chinen\footnotemark[1]}
\date{}
\maketitle

\begin{abstract}
The formal weight enumerators were first introduced by M. Ozeki. They form a ring of invariant polynomials which is similar to that of the weight enumerators of Type II codes. Later, the zeta functions for linear codes were discovered and their theory was developed by I. Duursma. It was generalized to certain invariant polynomials including Ozeki's formal weight enumerators by the present author. One of the famous and important problems is whether extremal weight enumerators satisfy the Riemann hypothesis. In this paper, first we formulate the notion of divisible formal weight enumerators and propose an algorithm for the efficient search of the formal weight enumerators divisible by two. The main tools are the binomial moments. It leads to the discovery of several new families of formal weight enumerators. Then, as a result, we find examples of extremal formal weight enumerators which do not satisfy the Riemann hypothesis. 
\end{abstract}%  words

\footnotetext[1]{Department of Mathematics, School of Science and Engineering, Kindai University. 3-4-1, Kowakae, Higashi-Osaka, 577-8502 Japan. E-mail: chinen@math.kindai.ac.jp}
%\footnotetext{This work was supported by JSPS KAKENHI Grant Number JP26400028. }

%\footnotetext[1]{}

%%%%%%%%%%%%%%%%%%%%%%%%%%%%%%%%%%%%%%%%%%%%%%%%%%%%%%%%%%%%%
\section{Introduction}\label{section:intro}
The notion of the formal weight enumerator was first introduced to number theory and coding theory by Ozeki \cite{Oz} in 1997: he called a polynomial $W(x,y)=\sum_{i=0}^n A_i x^{n-i}y^i \in {\bf C}[x,y]$ a formal weight enumerator if 
\begin{eqnarray}\label{eq:fwe_sigma2}
W^{\sigma_2}(x,y)&=&-W(x,y),\\
\label{eq:fwe_tau}
W^\tau(x,y)&=&W(x,y)
\end{eqnarray}
are satisfied, where 
$$\sigma_2=\frac{1}{\sqrt 2}\left(\begin{array}{rr} 1 & 1 \\ 1 & -1 \end{array}\right),\quad
\tau=\left(\begin{array}{rr} 1 & 0 \\ 0 & i \end{array}\right) 
\quad (i=\sqrt{-1})$$
and the action of a matrix $\sigma=\left(\begin{array}{cc} a & b \\ c & d \end{array}\right)$ on a polynomial $f(x,y)\in {\bf C}[x,y]$ is defined by $f^\sigma(x,y)=f(ax+by, cx+dy)$. Ozeki's formal weight enumerators are very close to the weight enumerators of so-called Type II codes, an important class of self-dual codes which are defined over the finite field ${\ff 2}$ (see the discussion that follows Theorem \ref{thm:gleason_pierce}). In fact, the weight enumerator $W_C(x,y)$ of a Type II code $C$ is characterized by 
\begin{eqnarray}\label{eq:type2_sigma2}
{W_C}^{\sigma_2}(x,y)&=&W_C(x,y),\\
\label{eq:type2_tau}
{W_C}^\tau(x,y)&=&W_C(x,y). 
\end{eqnarray}
The invariance by $\sigma_2$ means that $C$ is self-dual over ${\ff 2}$, the invariance by $\tau$ means that $C$ is ``divisible by four'' (the weights of all the codewords are divisible by four). The weight enumerators of Type II codes belong to the invariant polynomial ring
$$R_{\rm II} := {\bf C}[x,y]^{G_9}={\bf C}[W_{{\mathcal H}_{8}}(x,y), W_{{{\mathcal G}_{24}}}(x,y)],$$
where $G_9=\langle \sigma_2, \tau\rangle$ (the group No. 9 in Shephard-Todd \cite{ShTo}), 
$$W_{{\mathcal H}_{8}}(x,y)=x^8+14x^4y^4+y^8$$
(the weight enumerator of the extended Hamming code) and
$$W_{{{\mathcal G}_{24}}}(x,y)=x^{24}+759x^{16}y^8+2576x^{12}y^{12}+759x^8y^{16}+y^{24}$$
(that of the extended Golay code, see for example, Conway-Sloane \cite[p.192]{CoSl}). On the other hand, Ozeki's formal weight enumerators are contained in the invariant polynomial ring
$$R_{\rm II}^- :={\bf C}[x,y]^{G_8}={\bf C}[W_{{\mathcal H}_{8}}(x,y), W_{12}(x,y)],$$
where $G_8$ is the group No. 8 in Shephard-Todd \cite{ShTo}: 
$$G_8=\left\langle 
\frac{1-i}{2}\left(\begin{array}{rr} 1 & -1 \\ 1 & 1 \end{array}\right),\ 
\left(\begin{array}{rr} -i & 0 \\ 0 & 1 \end{array}\right) 
\right\rangle$$
and 
\begin{equation}\label{eq:W12}
W_{12}(x,y)=x^{12}-33x^8y^4-33x^4y^8+y^{12}.
\end{equation}
Ozeki \cite{Oz} used the polynomial $W_{12}(x,y)$ to construct the Eisenstein series $E_6(z)$ by the Brou\'e-Enguehard map. 

As to the divisibility of the self-dual codes, the following theorem is well-known (see Sloane \cite[Section 6]{Sl}):
%%%%%%%%%%%%%%%%%%%%%% theorem %%%%%%%%%%%%%%%%%%%%%%
\begin{thm}[Gleason-Pierce]\label{thm:gleason_pierce}
Suppose a self-dual code over $\ff q$ is divisible by $c>1$, that is, the weight of any codeword is divisible by $c$. Then, $(q,c)$ must be one of the following: 
$$(q,c)=(2,2), (2,4), (3,3), (4,2)$$
or $q$ is arbitrary and $c=2$. 
\end{thm}
A self-dual code over ${\ff 2}$ divisible two (the case $(q,c)=(2,2)$) is called a ``Type I code''. The other cases $(q,c)=(2,4), (3,3)$ and $(4,2)$ are called Types II, III and IV, respectively. In the last case, that is, the case where $q$ is arbitrary and $c=2$, the weight enumerator is given by powers of 
\begin{equation}\label{eq:deg2_inv}
W_{2,q}(x,y)=x^2+(q-1)y^2
\end{equation}
only, and this is rather a trivial case ($W_{2,q}(x,y)$ plays an important role in the present paper, though, see the remark before Corollary \ref{cor:moment_odd_deg}). So we know that (non-trivial) divisible self-dual codes exist only for very few pairs of $(q,c)$. 

However, removing the structure of linear codes and allowing $q$ to be any positive real number other than one, we will find many (possibly infinite) families of polynomials which resemble Ozeki's formal weight enumerators (for example, for $q=4/3, 4\pm 2\sqrt{2}, 2\pm 2\sqrt{5}/5, 8\pm 4\sqrt{3}$, etc., see Section \ref{section:examples}). So we begin by redefining the formal weight enumerator as follows: 
%%%%%%%%%%%%%%%%%%%%%% defintion %%%%%%%%%%%%%%%%%%%%%%
\begin{dfn}\label{dfn:fwe}
We call a homogeneous polynomial 
\begin{equation}\label{eq:fwe_form}
W(x,y)=x^n+\sum_{i=d}^n A_i x^{n-i} y^i\in {\bf C}[x,y]\quad (A_d\ne 0)
\end{equation}
a formal weight enumerator if 
\begin{equation}\label{eq:fwe_transf}
W^{\sigma_q}(x,y)=-W(x,y)
\end{equation}
for some $q\in {\bf R}$, $q> 0$, $q\ne1$, where
\begin{equation}\label{eq:macwilliams}
\sigma_q=\frac{1}{\sqrt{q}}\left(\begin{array}{rr} 1 & q-1 \\ 1 & -1 \end{array}\right).
\end{equation}
Moreover, for some fixed $c\in {\bf N}$, we call $W(x,y)$ divisible by $c$ if 
$$A_i\ne 0\quad \Rightarrow \quad c|i$$
is satisfied. 
\end{dfn}
The transformation which is defined by $\sigma_q$ is often called the MacWilliams transform. 

There are two aims in this paper. The first one is to propose an algorithm for the search of divisible formal weight enumerators $W(x,y)$. We restrict ourselves to the case where $c=2$. Our main tools are the binomial moments (MacWilliams-Sloane \cite[pp.130-131]{MaSl}). It is surprising that the possible values of $q$ are determined for each given degree of $W(x,y)$ by the identities satisfied by its binomial moments, due to the condition $W^{\sigma_q}(x,y)=-W(x,y)$. We can use this fact to establish an algorithm to look for the candidates of $q$. In connection to the Gleason-Pierce theorem, a necessary condition for existence of a self-dual code over $\ff q$ divisible by $c$ is given in Sloane \cite[(6.3.3)]{Sl}, that is, $(1-\alpha)^2/q$ is an algebraic integer, where $\alpha$ is a primitive $c$-th root of unity. It is a necessary condition also for existence of divisible formal weight enumerators in Definition \ref{dfn:fwe}, but is a little too weak to specify the values of $q$, for lack of the restriction that $q$ is a prime power. Our method in this paper can pick up the candidates of $q$ effectively. 

Our second aim is related to the theory of zeta functions for linear codes. The theory was discovered and developed by Duursma \cite{Du1} -- \cite{Du4}, later the author \cite{Ch1} and \cite{Ch2} generalized it to the cases of some invariant polynomial rings, including $R_{\rm II}^-$. Among the problems concerning zeta functions for linear codes, the following one is very famous and important (\cite[Open Problem 4.2]{Du3}): 
%%%%%%%%%%%%%%%%%%%%%% problem %%%%%%%%%%%%%%%%%%%%%%
\begin{prob}[Duursma]\label{prob:duursma_extremal}
Prove or disprove that all extremal weight enumerators satisfy the Riemann hypothesis. 
\end{prob}
For the definition of the extremal weight enumerator, see Definition \ref{dfn:extremal_code}, and for that of the Riemann hypothesis, see Definition \ref{dfn:RH}. As far as the divisible self-dual codes (see Theorem \ref{thm:gleason_pierce}) are concerned, we do not know any example of an extremal weight enumerator not satisfying the Riemann hypothesis. However, in the course of our search for divisible formal weight enumerators, we encounter some examples which are extremal but do {\it not} satisfy the Riemann hypothesis. These examples may suggest that it is not only the extremal property that has an influence on the truth of the Riemann hypothesis. Our second aim in this paper is to pose a question about the relation between the extremal property and the Riemann hypothesis by showing such examples (see Section \ref{section:extremal_no_rh}). 

The algorithm which is proposed in this paper can be applied to the search of polynomials of odd degrees, invariant by $\sigma_q$ and divisible by two. It leads to finding another polynomial ring, in which we can observe similar phenomena, that is, the discovery of extremal invariant polynomials not satisfying the Riemann hypothesis (the reader is referred to \cite{Ch3}). 

The rest of the paper is organized as follows: in Section \ref{section:moments}, we introduce the identities satisfied by the binomial moments of formal weight enumerators divisible by two and propose an algorithm for the search of them. In Section \ref{section:examples}, we calculate some examples of formal weight enumerators. Section \ref{section:zeta} is devoted to a brief summary of the zeta functions for linear codes and invariant polynomials, and we give some new examples of formal weight enumerators satisfying the Riemann hypothesis. In Section \ref{section:extremal_no_rh}, we give several examples of extremal formal weight enumerators not satisfying the Riemann hypothesis. In the last section, we give some remarks and problems. 
%%%%%%%%%%%%%%%%%%%%%%%%%%%%%%%%%%%%%%%%%%%%%%%%%%%%
%%%%%%%%%%%%%%%%%%%%%%%%%%%%%%%%%%%%%%%%%%%%%%%%%%%%
\section{Binomial moments for formal weight enumerators}\label{section:moments}
In this section, we propose an algorithm for the search of divisible formal weight enumerators. Our starting point is the following: 
%%%%%%%%%%%%%%%%%% theorem %%%%%%%%%%%%%%%%%%
\begin{thm}[Binomial moments]\label{thm:moment}
Let
\begin{eqnarray*}
W(x,y)&=&\sum_{i=0}^n A_i x^{n-i} y^i\\
&=&x^n+\sum_{i=d}^n A_i x^{n-i} y^i \qquad (A_d\ne 0) 
\end{eqnarray*}
be a formal weight enumerator satisfying $W^{\sigma_q}(x,y)=-W(x,y)$. Then we have
\begin{equation}\label{eq:moment_fwe}
\sum_{i=0}^{n-\nu} {{n-i}\choose{\nu}} A_i = -q^{\frac{n}{2}-\nu} \sum_{i=0}^\nu {{n-i}\choose{n-\nu}} A_i
\qquad(\nu=0,1,\cdots, n).
\end{equation}
\end{thm}
\prf We have
\begin{equation}\label{eq:macwilliams_fwe}
\sum_{i=0}^n A_i x^{n-i} y^i = -\frac{1}{q^{n/2}}\sum_{i=0}^n A_i (x+(q-1)y)^{n-i}(x-y)^i
\end{equation}
by $W(x,y)=-W^{\sigma_q}(x,y)$. Let $y=1$ in (\ref{eq:macwilliams_fwe}) and differentiate it $\nu$ times with respect to $x$, next let $x=1$. Then we get (\ref{eq:moment_fwe}) (see also MacWilliams-Sloane \cite[p.131, Problem (6)]{MaSl}). \qed

\medskip
\noindent We divide our search into two cases according to the parity of $\deg W(x,y)$. If $\deg W(x,y)$ is even, we assume 
\begin{equation}\label{eq:fwe_even_deg}
W(x,y)=\sum_{i=0}^n A_i x^{2(n-i)} y^{2i}\qquad(A_0=1).
\end{equation}
For odd degrees, we assume 
\begin{equation}\label{eq:fwe_odd_deg}
W(x,y)=\sum_{i=0}^n A_i x^{2(n-i)+1} y^{2i}\qquad(A_0=1).
\end{equation}
(I) {\it The even degree case}

\medskip
\noindent The identities satisfied by the binomial moments of (\ref{eq:fwe_even_deg}) are easily obtained from Theorem \ref{thm:moment}: 
%%%%%%%%%%%%%%%%%% corollary %%%%%%%%%%%%%%%%%%
\begin{cor}\label{cor:moment_even_deg}
Suppose the polynomial (\ref{eq:fwe_even_deg}) satisfies $W^{\sigma_q}(x,y)=-W(x,y)$. Then we have 
\begin{equation}\label{eq:moment_even_deg}
\sum_{i=0}^{n-[\frac{\nu+1}{2}]} {{2n-2i}\choose{\nu}} A_i = 
-q^{n-\nu} \sum_{i=0}^{[\frac{\nu}{2}]} {{2n-2i}\choose{2n-\nu}} A_i
\qquad(\nu=0,1,\cdots, 2n),
\end{equation}
where $[x]$ means the greatest integer not exceeding x. 
\end{cor}

\medskip
\noindent The formula (\ref{eq:moment_even_deg}) gives $2n+1$ linear equations of $A_0, A_1, \cdots, A_n$, but the cases $\nu=n+1, n+2, \cdots, 2n$ are essentially the same as the cases $\nu=n-1, n-2, \cdots, 0$, respectively. So it suffices to consider the cases $\nu=0,1, \cdots, n$. Note that the formula (\ref{eq:moment_even_deg}) remains nontrivial when $\nu=n$, due to the condition $W^{\sigma_q}(x,y)=-W(x,y)$. Indeed, it becomes
$$2\sum_{i=0}^{[\frac{n}{2}]} {{2n-2i}\choose{n}} A_i = 0.$$
This is the striking difference from the case of self-dual weight enumerators (we have  ${W_C}^{\sigma_q}(x,y)=W_C(x,y)$ and $\deg W_C(x,y)$ is automatically even), in which case the formula corresponding to (\ref{eq:moment_even_deg}) with $\nu=n$ disappears since $\sum_{i=0}^{[n/2]} {{2n-2i}\choose{n}} A_i = \sum_{i=0}^{[n/2]} {{2n-2i}\choose{n}} A_i$. 

Thus we get a system of $n+1$ homogeneous linear equations of $n+1$ unknowns $A_0, A_1, \cdots, A_n$. Let $A(n,q)$ be the coefficient matrix of this system. 

\medskip
\noindent{\bf An Algorithm}

\nopagebreak
\medskip
\noindent We can use the above system of linear equations to look for a formal weight enumerator of degree $2n$, divisible by two. Since we need $A_0=1$, the system must have a nontrivial solution. So first we find a number $q$ satisfying $|A(n,q)|=0$. Next we find other coefficients $A_1, A_2, \cdots, A_n$ and construct  $W(x,y)=\sum_{i=0}^n A_i x^{2(n-i)} y^{2i}$. Finally, we verify $W^{\sigma_q}(x,y)=-W(x,y)$. 

\medskip
\rem As we will see later, once we get a formal weight enumerator $W(x,y)$ for some $q=q_0$ and a certain $n=n_0$, then we have $|A(n, q_0)|=0$ for all $n\geq n_0$. This is due to the existence of the polynomial $W_{2,q_0}(x,y)$ of (\ref{eq:deg2_inv}). In fact, if deg$W(x,y)=2n$, then $W(x,y)W_{2,q_0}(x,y)^k$ ($k\geq 0$) is a formal weight enumerator of degree $2n+2k$ for $q=q_0$. 

\medskip
\noindent(II) {\it The odd degree case}

\medskip
This case is similar. The identities satisfied by the binomial moments of (\ref{eq:fwe_odd_deg}) are the following: 
%%%%%%%%%%%%%%%%%% corollary %%%%%%%%%%%%%%%%%%
\begin{cor}\label{cor:moment_odd_deg}
Suppose the polynomial (\ref{eq:fwe_odd_deg}) satisfies $W^{\sigma_q}(x,y)=-W(x,y)$. Then we have 
\begin{equation}\label{eq:moment_odd_deg}
\sum_{i=0}^{[\frac{2n+1-\nu}{2}]} {{2n+1-2i}\choose{\nu}} A_i = 
-q^{n-\nu+1/2} \sum_{i=0}^{[\frac{\nu}{2}]} {{2n+1-2i}\choose{2n+1-\nu}} A_i
\qquad(\nu=0,1,\cdots, 2n+1),
\end{equation}
where $[x]$ means the greatest integer not exceeding x. 
\end{cor}
The formula (\ref{eq:moment_odd_deg}) gives essentially $(2n+2)/2=n+1$ linear equations of $n+1$ unknowns $A_0, A_1, \cdots, A_n$. So we can apply the same reasoning as the even degree case. We denote the coefficient matrix of the system of the linear equations by $B(n,q)$. 
%%%%%%%%%%%%%%%%%%%%%%%%%%%%%%%%%%%%%%%%%%%%%%%%%%%%
%%%%%%%%%%%%%%%%%%%%%%%%%%%%%%%%%%%%%%%%%%%%%%%%%%%%
\section{Some examples}\label{section:examples}
We follow the classification in the last section. 

\medskip
\noindent(I) {\it The even degree case}

\nopagebreak
\medskip Recall that $A(n,q)$ is the coefficient matrix of the system of $n+1$ linear equations 
$$\sum_{i=0}^{n-[\frac{\nu+1}{2}]} {{2n-2i}\choose{\nu}} A_i + q^{n-\nu} \sum_{i=0}^{[\frac{\nu}{2}]} {{2n-2i}\choose{2n-\nu}} A_i =0
\qquad(\nu=0,1,\cdots, n)$$
of $n+1$ unknowns $A_0, A_1, \cdots, A_n$. 

\medskip
\noindent (i) {\it The case $n=1$}

\medskip
\noindent We can easily see that 
$$A(1,q)=\left(\begin{array}{cc} 1+q & 1 \\ 4 & 0 \end{array}\right),$$
which is regular for all $q$, so there is no formal weight enumerator of degree two divisible by two. 

\medskip
\noindent (ii) {\it The case $n=2$}

\medskip
\noindent We can see 
$$A(2,q)=\left(\begin{array}{ccc} 1+q^2 & 1 & 1 \\ 4(1+q) & 2 & 0 \\ 12 & 2 & 0\end{array}\right)$$
and $|A(2,q)|=8(q-2)$. So we have $q=2$. Setting $A_0=1$, we have other coefficients $A_1=-6$ and $A_2=1$. Let 
\begin{equation}\label{eq:fwe_4_2}
\varphi_4(x,y)=x^4-6x^2y^2+y^4.
\end{equation}
We can verify that ${\varphi_4}^{\sigma_2}(x,y)=-\varphi_4(x,y)$, thus we have found a non-trivial formal weight enumerator for $(q,c)=(2,2)$. The ring 
$$R_{\rm I}^-:={\bf C}[W_{2,2}(x,y)=x^2+y^2, \varphi_4(x,y)]$$
is, so to speak, the ring of ``Type I formal weight enumerators'' in comparison with $R_{\rm I}$, the ring of the weight enumerators of Type I codes (the self-dual codes over ${\ff 2}$ divisible by two): 
$$R_{\rm I}={\bf C}[W_{2,2}(x,y), W_{{\mathcal H}_{8}}(x,y)]$$
(see for example, Conway-Sloane \cite[p.186]{CoSl}). Obviously, $R_{\rm I}^-$ is the largest ring that contains the polynomials $W(x,y)$ divisible by two with the property $W^{\sigma_2}(x,y)=\pm W(x,y)$. It contains various well known polynomials, for example, 
\begin{eqnarray}
W_{{\mathcal H}_{8}}(x,y)&=&\frac{1}{4}\left(3W_{2,2}(x,y)^4+\varphi_4(x,y)^2 \right),\nonumber\\
W_{12}(x,y)&=&\frac{1}{8}\left(9W_{2,2}(x,y)^4 \varphi_4(x,y)-\varphi_4(x,y)^3 \right),\label{eq:12lin_comb}\\
W_{{{\mathcal G}_{24}}}(x,y)&=&\frac{1}{128}\left(33W_{2,2}(x,y)^{12}+96W_{2,2}(x,y)^8 \varphi_4(x,y)^2-3W_{2,2}(x,y)^4 \varphi_4(x,y)^4 \right.\nonumber\\
& & \left. +2\varphi_4(x,y)^6\right).\nonumber
\end{eqnarray}
The formal weight enumerators in $R_{\rm I}^-$ are investigated more closely in our subsequent paper \cite{Ch4}. 

\medskip
\noindent (iii) {\it The case $n=3$}

\medskip
\noindent We have
$$A(3,q)=\left(\begin{array}{cccc}
1+q^3 & 1 & 1 & 1 \\
6(1+q^2) & 4 & 2 & 0 \\
15(1+q) & 6+q & 1 & 0 \\
40 & 8 & 0 & 0 
\end{array}\right)$$
and $|A(3,q)|=16(q-2)(3q-4)$. Other than $q=2$, we get $q=4/3$. Setting $A_0=1$, we obtain $A_1=-5$, $A_2=5/3$, $A_3=-1/27$. We can verify (with some computer algebra system) that 
$$\varphi_6(x,y)=x^6-5x^4y^2+\frac{5}{3}x^2y^4-\frac{1}{27}y^6$$
satisfies ${\varphi_6}^{\sigma_{4/3}}(x,y)=-\varphi_6(x,y)$. 

\medskip
\noindent (iv) {\it The case $n=4$}

\medskip
\noindent We have
$$A(4,q)=\left(\begin{array}{ccccc}
1+q^4 & 1 & 1 & 1 & 1 \\
8(1+q^3) & 6 & 4 & 2 & 0 \\
28(1+q^2) & 15+q^2 & 6 & 1 & 0 \\
56(1+q) & 20+6q & 4 & 0 & 0 \\
140 & 30 & 2 & 0 & 0 
\end{array}\right)$$
and $|A(4,q)|=32(q-2)(3q-4)(q^2-8q+8)$. We newly find $q=4\pm 2\sqrt{2}$. Setting $A_0=1$, we obtain 
\begin{eqnarray}
\varphi_8^{\pm}(x,y) &=& x^8-(84\pm 56\sqrt{2})x^6y^2+(1190\pm 840\sqrt{2})x^4y^4\nonumber\\
& & -(2772\pm 1960\sqrt{2})x^2y^6+(577\pm 408\sqrt{2})y^8\label{eq:fwe_deg8}
\end{eqnarray}
with the property ${\varphi_8^{\pm}}^{\sigma_{4\pm 2\sqrt{2}}}(x,y)=-\varphi_8^{\pm}(x,y)$. 

\medskip
\noindent (v) {\it The case $n=5$}

\nopagebreak
\medskip
\noindent We omit the explicit form of the matrix $A(5,q)$. For the determinant, we have 
$$|A(5,q)|=-64(q-2)(3q-4)(q^2-8q+8)(5q^2-20q+16).$$
The new values of $q$ are $q=2 \pm 2\sqrt{5}/5$. We can verify that the polynomials 
\begin{eqnarray}
\varphi_{10}^{\pm}(x,y) &=& x^{10}-(45\pm 18\sqrt{5})x^8y^2+(378\pm 168\sqrt{5})x^6y^4 -\left(714 \pm \frac{1596}{5}\sqrt{5}\right)x^4y^6\nonumber\\
& &  +\left( \frac{1449}{5} \pm \frac{648}{5}\sqrt{5} \right)x^2y^8 -\left( \frac{61}{5} \pm \frac{682}{125}\sqrt{5} \right) y^{10} \label{eq:fwe_deg10}
\end{eqnarray}
satisfy ${\varphi_{10}^{\pm}}^{\sigma_{2\pm 2\sqrt{5}/5}}(x,y)=-\varphi_{10}^{\pm}(x,y)$. 

\medskip
\noindent (vi) {\it The case $n=6$}

\medskip
\noindent We have 
$$|A(6,q)|=128(q-2)^2(3q-4)(q^2-8q+8)(5q^2-20q+16)(q^2-16q+16).$$
The new values of $q$ are $q=8\pm 4\sqrt{3}$. We find the following formal weight enumerators of degree 12 for $\sigma_{8\pm 4\sqrt{3}}$: 
\begin{eqnarray}
\varphi_{12}^{\pm}(x,y) &=& x^{12}-(462 \pm 264\sqrt{3})x^{10}y^2+(48015 \pm 27720\sqrt{3})x^8y^4\nonumber\\
& & -(1248324 \pm 720720\sqrt{3})x^6y^6+(9314415 \pm 5377680\sqrt{3})x^4y^8\nonumber\\
& & -(17297742 \pm 9986856\sqrt{3})x^2y^{10}+(3650401 \pm 2107560\sqrt{3})y^{12}. \label{eq:fwe_deg12}
\end{eqnarray}

\medskip
\noindent(II) {\it The odd degree case}

\nopagebreak
\medskip Recall that $B(n,q)$ is the coefficient matrix of the system of $n+1$ linear equations 
$$\sum_{i=0}^{[\frac{2n+1-\nu}{2}]} {{2n+1-2i}\choose{\nu}} A_i + q^{n-\nu+1/2} \sum_{i=0}^{[\frac{\nu}{2}]} {{2n+1-2i}\choose{2n+1-\nu}} A_i =0
\qquad(\nu=0,1,\cdots, n)$$
of $n+1$ unknowns $A_0, A_1, \cdots, A_n$. 

\medskip
\noindent (i) {\it The case $n=1$}

\medskip
\noindent We have 
$$B(1,q)=\left(\begin{array}{cc} 1+q\sqrt{q} & 1 \\ 3(1+\sqrt{q}) & 1\end{array}\right)$$
and $|B(1,q)|=q\sqrt{q}-3\sqrt{q}-2$. So we get $q=4$ from $|B(1,q)|=0$. Setting $A_0=1$, we obtain $A_1=-9$.  Let 
\begin{equation}\label{eq:fwe_3_4}
\varphi_3(x,y)=x^3-9xy^2.
\end{equation}
Then we can easily verify ${\varphi_3}^{\sigma_4}(x,y)=-\varphi_3(x,y)$, so we have found a non-trivial formal weight enumerator for $(q,c)=(4,2)$. Similarly to the case of $R_{\rm I}$ and $R_{\rm I}^-$, we can construct the ring of ``Type IV formal weight enumerators''
$$R_{\rm IV}^- :={\bf C}[W_{2,4}(x,y)=x^2+3y^2, \varphi_3(x,y)].$$
The ring of Type IV weight enumerators is 
$$R_{\rm IV} :={\bf C}[W_{2,4}(x,y), x^6+45x^2y^4+18y^6]$$
(see for example, Conway-Sloane \cite[p.203]{CoSl}). The formal weight enumerators in $R_{\rm IV}^-$ are investigated more closely in our subsequent paper \cite{Ch4}. 

\medskip
\noindent (ii) {\it The case $n=2$}

\medskip
\noindent We have 
$$B(2,q)=
\left(\begin{array}{ccc} 1+q^2\sqrt{q} & 1 & 1 \\
                        5(1+q\sqrt{q}) & 3 & 1 \\
                        10(1+\sqrt{q}) & 3+\sqrt{q} & 0\end{array}\right)$$
and putting $\sqrt{q}=t$, we have 
$$|B(2,q)|=-(t+1)^3(t-2)(t^2+2t-4).$$
From $t>0$, we newly find $t=-1+\sqrt{5}$ and $q=t^2=6-2\sqrt{5}$. Setting $A_0=1$, we obtain $A_1=-50+20\sqrt{5}$, $A_2=225-100\sqrt{5}$. We can verify that 
\begin{equation}\label{eq:fwe_deg5}
\varphi_5(x,y)=x^5+(-50+20\sqrt{5})x^3y^2+(225-100\sqrt{5})xy^4
\end{equation}
satisfies ${\varphi_5}^{\sigma_{6-2\sqrt{5}}}(x,y)=-{\varphi_5}(x,y)$. This is a non-trivial formal weight enumerator for $(q,c)=(6-2\sqrt{5}, 2)$. 

%%%%%%%%%%%%%%%%%%%%%%%%%%%%%%%%%%%%%%%%%%%%%%%%%%%%
%%%%%%%%%%%%%%%%%%%%%%%%%%%%%%%%%%%%%%%%%%%%%%%%%%%%
\section{Zeta functions and the Riemann hypothesis for formal weight enumerators}\label{section:zeta}
In this section, we summarize some basic definitions and facts on zeta functions for codes and for formal weight enumerators. We always assume $d, d^\perp \geq 2$ where $d^\perp$ is defined by 
$$W^{\sigma_q}(x,y)=\pm x^n + A_{d^\perp} x^{n-d^\perp} y^{d^\perp}+ \cdots,$$
when considering the zeta functions for $W(x,y)$ of the form (\ref{eq:fwe_form}) (see \cite[Section 2]{Du2}). 
%%%%%%%%%%%%%%%%%% definition %%%%%%%%%%%%%%%%%%
\begin{dfn}\label{dfn:zeta}
For any homogeneous polynomial of the form (\ref{eq:fwe_form}) and $q\in{\bf R}$ ($q>0, q\ne 1$), there exists a unique polynomial $P(T)\in{\bf C}[T]$ of degree at most $n-d$ such that
\begin{equation}\label{eq:zeta_duursma}
\frac{P(T)}{(1-T)(1-qT)}(y(1-T)+xT)^n=\cdots +\frac{W(x,y)-x^n}{q-1}T^{n-d}+ \cdots.
\end{equation}
We call $P(T)$ and $Z(T)=P(T)/(1-T)(1-qT)$ the zeta polynomial and the zeta function of $W(x,y)$, respectively. 
\end{dfn}
Zeta functions of this type was first defined by Duursma \cite{Du1} for $W(x,y)=W_C(x,y)$, the weight enumerator of a linear code $C$. In that case, the number $q$ must be chosen as the number of elements of the finite field over which $C$ is defined. But the definition can be easily extended to any polynomial of the form (\ref{eq:fwe_form}). For an elementary proof of existence of $P(T)$, see for example \cite[Appendix A]{Ch2}. 

If $W(x,y)$ is a formal weight enumerator of the form (\ref{eq:fwe_form}) with $W^{\sigma_q}(x,y)=-W(x,y)$, we can prove the following in a similar way to \cite[p.59]{Du2} (see also \cite[Theorem 2.1]{Ch1}): 
%%%%%%%%%%%%%%%%%% theorem %%%%%%%%%%%%%%%%%%
\begin{thm}[Functional equation]\label{thm:func_eq}
The zeta polynomial $P(T)$ of a formal weight enumerator $W(x,y)$ with $W^{\sigma_q}(x,y)=-W(x,y)$ is of degree $2g$ ($g=n/2+1-d$) and satisfies
\begin{equation}\label{eq:func_eq}
P(T)=-P\left(\frac{1}{qT}\right)q^g T^{2g}.
\end{equation}
\end{thm}
It is interesting to compare this with the case of the zeta polynomial $P_C(T)$ for a self-dual code $C$ over ${\ff q}$: 
$$P_C(T)=P_C\left(\frac{1}{qT}\right)q^g T^{2g}$$
(see \cite[Section 2]{Du2} and \cite[Section 4]{Du3}). In the case of codes, $g$ is called the genus of $C$. For a formal weight enumerator $W(x,y)$, we also call $g$ the genus of $W(x,y)$. Note that we have
\begin{equation}\label{eq:ineq_genus}
d\leq \frac{n}{2}+1
\end{equation}
because $g$ must satisfy $g\geq 0$ (otherwise, $P(T)$ would not be a polynomial). 

Now we can formulate the Riemann hypothesis for (formal) weight enumerators (see also \cite[Section 2]{Ch1}): 
%%%%%%%%%%%%%%%%%% definition %%%%%%%%%%%%%%%%%%
\begin{dfn}[Riemann hypothesis]\label{dfn:RH}
A (formal) weight enumerator $W(x,y)$ with $W^{\sigma_q}(x,y)=\pm W(x,y)$ satisfies the Riemann hypothesis if all the zeros of $P(T)$ have the same absolute value $1/\sqrt{q}$. 
\end{dfn}

In the rest of this section, we show some examples of formal weight enumerators that satisfy the Riemann hypothesis, among those were found in the last section. 

\medskip
\noindent(i) {\it The ring $R_{\rm I}^-$}

\medskip
\noindent The ring $R_{\rm I}^-$ is defined by 
$$R_{\rm I}^-:={\bf C}[W_{2,2}(x,y)=x^2+y^2, \varphi_4(x,y)=x^4-6x^2y^2+y^4].$$ 
The formal weight enumerators in $R_{\rm I}^-$ are of the forms
$$W_{2,2}(x,y)^l \varphi_4(x,y)^{2m+1}\qquad (l,m\geq 0)$$
and their suitable linear combinations (an example of such a linear combination is given in (\ref{eq:12lin_comb})). We show some pairs of formal weight enumerators $W(x,y)$ and their zeta polynomials $P(T)$. We can easily verify that all the roots of $P(T)$ lie on the circle $|T|=1/\sqrt{2}$ in each case: 
\begin{eqnarray*}
W(x,y) &=& \varphi_4(x,y)=x^4-6x^2y^2+y^4, \\
P(T) &=& 2T^2-1;\\
W(x,y) &=& W_{2,2}(x,y) \varphi_4(x,y)=x^6-5x^4y^2-5x^2y^4+y^6, \\
P(T) &=& \frac{1}{3}(2T^2-1)(2T^2+1); \\
W(x,y) &=& W_{2,2}(x,y)^2 \varphi_4(x,y)=x^8-4x^6y^2-10x^4y^4-4x^2y^6+y^8,\\
P(T) &=& \frac{1}{7}(2T^2-1)(4T^4+2T^2+1);\\
W(x,y) &=& W_{2,2}(x,y)^3 \varphi_4(x,y)=x^{10}-3x^8y^2-14x^6y^4-14x^4y^6-3x^2y^8+y^{10},\\
P(T) &=& \frac{1}{15}(2T^2-1)(2T^2+1)(2T^2+2T+1)(2T^2-2T+1);\\
W(x,y) &=& W_{12}(x,y)=x^{12}-33x^8y^4-33x^4y^8+y^{12}\quad (\mbox{see }(\ref{eq:12lin_comb})),\\
P(T) &=& \frac{1}{15}(2T^2-1)(2T^2+1)(2T^2+2T+1).
\end{eqnarray*}

\medskip
\noindent(ii) {\it The ring $R_{\rm IV}^-$}

\medskip
\noindent The ring $R_{\rm IV}^-$ is defined by 
$$R_{\rm IV}^- :={\bf C}[W_{2,4}(x,y)=x^2+3y^2, \varphi_3(x,y)=x^3-9xy^2].$$
In the following, we can easily verify that all the roots of $P(T)$ lie on the circle $|T|=1/2$ in each case: 
\begin{eqnarray*}
W(x,y) &=& \varphi_3(x,y)=x^3-9xy^2, \\
P(T) &=& 2T-1;\\
W(x,y) &=& W_{2,4}(x,y)\varphi_3(x,y)=x^5-6x^3y^2-27xy^4, \\
P(T) &=& \frac{1}{5}(2T-1)(4T^2+1); \\
W(x,y) &=& W_{2,4}(x,y)^2\varphi_3(x,y)=x^7-3x^5y^2-45x^3y^4-81xy^6,\\
P(T) &=& \frac{1}{21}(2T-1)(4T^2-2T+1)(4T^2+2T+1);\\
W(x,y) &=& W_{2,4}(x,y)^3\varphi_3(x,y)=x^9-54x^5y^4-216x^3y^6-243xy^8,\\
P(T) &=& \frac{1}{7}(2T-1)(4T^2+2T+1);\\
W(x,y) &=& \frac{1}{9}(8W_{2,4}(x,y)^4\varphi_3(x,y)+W_{2,4}(x,y)\varphi_3(x,y)^3)\\
       & & =x^{11}-30x^7y^4-336x^5y^6-1035x^3y^8-648xy^{10},\\
P(T) &=& \frac{1}{33}(2T-1)(16T^4+8T^3+6T^2+2T+1).
\end{eqnarray*}

\medskip
\noindent(iii) {\it Polynomials for $q=4/3$}

\medskip
\noindent We consider the ring $R_{4/3}^- :={\bf C}[W_{2,4/3}(x,y), \varphi_6(x,y)]$ where 
\begin{eqnarray*}
W_{2,4/3}(x,y) &=& x^2+\frac{1}{3}y^2,\\
\varphi_6(x,y) &=& x^6-5x^4y^2+\frac{5}{3}x^2y^4-\frac{1}{27}y^6.
\end{eqnarray*}
The formal weight enumerators in $R_{4/3}^-$ are of the forms
$$W_{2,4/3}(x,y)^l \varphi_6(x,y)^{2m+1}\qquad (l,m\geq 0)$$
and their suitable linear combinations. We can also consider the subring 
$$R_{4/3}={\bf C}[W_{2,4/3}(x,y), \varphi_6(x,y)^2]$$
to which the polynomials invariant under $\sigma_{4/3}$ belong. Invariant polynomials of the type of (\ref{eq:fwe_form}) are of the forms 
$$W_{2,4/3}(x,y)^l \varphi_6(x,y)^{2m}\qquad (l,m\geq 0, \ (l,m)\ne(0,0))$$
and their suitable linear combinations. In the following, we can easily verify that all the roots of $P(T)$ lie on the circle $|T|=\sqrt{3/4}=\sqrt{3}/2$ in each case: 
\begin{eqnarray*}
W(x,y) &=& \varphi_6(x,y) = x^6-5x^4y^2+\frac{5}{3}x^2y^4-\frac{1}{27}y^6, \\
P(T) &=& \frac{1}{9}(4T^2-3)(4T^2+2T+3);\\
W(x,y) &=& W_{2,4/3}(x,y)\varphi_6(x,y) = x^8-\frac{14}{3}x^6y^2+\frac{14}{27}x^2y^6-\frac{1}{81}y^8, \\
P(T) &=& \frac{1}{54}(4T^2-3)(16T^4+8T^3+15T^2+6T+9); \\
W(x,y) &=& W_{2,4/3}(x,y)^2 = x^4+\frac{2}{3}x^2y^2+\frac{1}{9}y^4\quad(\mbox{invariant under }\sigma_{4/3}),\\
P(T) &=& \frac{1}{9}(4T^2+2T+3).
\end{eqnarray*}
The formal weight enumerators and invariant polynomials in $R_{4/3}^-$ are investigated more closely in our subsequent paper \cite{Ch4}. 

\medskip
\noindent(iv) {\it Other polynomials}

\medskip
\noindent Among other polynomials, we only show the zeta polynomial for $\varphi_5(x,y)$ (see (\ref{eq:fwe_deg5})): 
\begin{eqnarray*}
W(x,y) &=& \varphi_5(x,y), \\
P(T) &=& \frac{2-\sqrt{5}}{4}(4T-(1+\sqrt{5}))(8T^2+4\sqrt{5}T+3+\sqrt{5}).
\end{eqnarray*}
All the roots of $P(T)$ are on the circle $|T|=1/\sqrt{q}=(1+\sqrt{5})/4$. 
See the next section for the Riemann hypothesis of formal weight enumerators for $q=4 \pm 2\sqrt{2}, 2 \pm 2\sqrt{5}/5$ and $8 \pm 4\sqrt{3}$. 
%%%%%%%%%%%%%%%%%%%%%%%%%%%%%%%%%%%%%%%%%%%%%%%%%%%%
%%%%%%%%%%%%%%%%%%%%%%%%%%%%%%%%%%%%%%%%%%%%%%%%%%%%
\section{Extremal formal weight enumerators not satisfying the Riemann hypothesis}\label{section:extremal_no_rh}
First we introduce the notion of extremal weight enumerators. For the first four cases of the Gleason-Pierce Theorem (Theorem \ref{thm:gleason_pierce}), that is the self-dual codes of Types I through IV, there are well-known upper bounds of the minimum distance $d$ from above by the code length $n$: 
%%%%%%%%%%%%%%%%%% theorem %%%%%%%%%%%%%%%%%%
\begin{thm}[Mallows-Sloane bound]\label{thm:Mallows-Sloane}
We have the following upper bounds for the minimum distance $d$ by the code length $n$: 
\begin{eqnarray*}
(\mbox{\rm Type I}) & & d\leq 2\left[\frac{n}{8}\right]+2, \\
(\mbox{\rm Type II}) & & d\leq 4\left[\frac{n}{24}\right]+4, \\
(\mbox{\rm Type III}) & & d\leq 3\left[\frac{n}{12}\right]+3, \\
(\mbox{\rm Type IV}) & & d\leq 2\left[\frac{n}{6}\right]+2.
\end{eqnarray*}
\end{thm}
\prf See for example, \cite[Theorem 3]{Du4} (see also \cite[p.138, Corollary]{Pl} and \cite[Chapter 19, Theorem 13]{MaSl}). \qed

\medskip
\noindent \rem For the Type I codes, a new (sharper) bound is established in Rains \cite{Ra}, but we adopt the classical bound since we are interested in the polynomials themselves. 

\medskip
\noindent We can define the extremal code as follows: 
%%%%%%%%%%%%%%%%%% definition %%%%%%%%%%%%%%%%%%
\begin{dfn}\label{dfn:extremal_code}
Among the codes of Types I through IV, the codes which attain the equality in the Mallows-Sloane bounds are called extremal codes. Weight enumerators of extremal codes are called extremal weight enumerators. 
\end{dfn}
From the practical point of view, extremal codes are good codes because the minimum distance is the largest among the given length $n$. The extremal weight enumerator is defined only in the words of polynomials of the form (\ref{eq:fwe_form}) without the structure of codes, so we can extend the notion to the case of formal weight enumerators (see also \cite[p.139]{Pl}): 
%%%%%%%%%%%%%%%%%% definition %%%%%%%%%%%%%%%%%%
\begin{dfn}\label{dfn:extremal_fwe}
Let 
$$S_{q,c,n}=\left\{W(x,y)\in {\bf C}[x,y]\ ; \ \begin{array}{l}
W(x,y)\mbox{ is of the form (\ref{eq:fwe_form})}, \\
W^{\sigma_q}(x,y)=-W(x,y),\\
W(x,y)\mbox{ is divisible by $c$}
\end{array} \right\}.$$
Then we call $W(x,y)\in S_{q,c,n}$ with maximal $d$ an extremal formal weight enumerator of degree $n$. 
\end{dfn}
In the case of formal weight enumerators, we can sometimes establish a similar bound to Theorem \ref{thm:Mallows-Sloane}, for example, 
%%%%%%%%%%%%%%%%%% theorem %%%%%%%%%%%%%%%%%%
\begin{thm}\label{thm:Mallows-Sloane_Ozeki_fwe}
Any formal weight enumerator of the form (\ref{eq:fwe_form}) in the ring $R_{\rm II}^-$ satisfies
\begin{equation}\label{eq:Mallows-Sloane_Ozeki_fwe}
d\leq 4 \left[\frac{n-12}{24}\right] + 4.
\end{equation}
\end{thm}
\prf See \cite[Theorem 3.4]{Ch1}. \qed

\medskip
\noindent Thus an extremal formal weight enumerator in $R_{\rm II}^-$ is a polynomial of the form (\ref{eq:fwe_form}) which satisfies $d=4[(n-12)/24]+4$. One of the examples is $W_{12}(x,y)$ in (\ref{eq:W12}). 

Even if we do not have a bound like (\ref{eq:Mallows-Sloane_Ozeki_fwe}) which is valid for all degrees $n$, we can find the extremal formal weight enumerators by simple manipulations of polynomials, once the degree $n$ is specified explicitly, as we will see in the discussion that follows. First we consider the cases $q=4\pm 2\sqrt{2}$. Let
\begin{eqnarray}
R_{4 + 2\sqrt{2}}^- &=& {\bf C}[W_{2, 4 + 2\sqrt{2}}(x,y), \varphi_8^+(x,y)], \label{eq:ring_fwe_4+2sqrt2}\\
R_{4 - 2\sqrt{2}}^- &=& {\bf C}[W_{2, 4 - 2\sqrt{2}}(x,y), \varphi_8^-(x,y)], \label{eq:ring_fwe_4-2sqrt2}
\end{eqnarray}
where $W_{2, 4 \pm 2\sqrt{2}}(x,y)=x^2+(3 \pm 2\sqrt{2})y^2$ (see (\ref{eq:fwe_deg8}) for the definitions of $\varphi_8^\pm (x,y)$). The rings $R_{4 + 2\sqrt{2}}^-$ and $R_{4 - 2\sqrt{2}}^-$ are those of formal weight enumerators for $q=4 + 2\sqrt{2}$ and $4 - 2\sqrt{2}$, respectively. Formal weight enumerators in $R_{4 + 2\sqrt{2}}^-$ are of the form
\begin{equation}\label{eq:form_fwe_4+2sqrt2}
W_{2, 4 + 2\sqrt{2}}(x,y)^l \varphi_8^+(x,y)^{2m+1} \qquad (l,m \geq 0)
\end{equation}
and their suitable linear combinations (the case of $q=4 - 2\sqrt{2}$ is similar). In the case of degree eight, the only possible pair of $(l,m)$ is $(l,m)=(0,0)$, so $\varphi_8^\pm (x,y)$ themselves are the extremal formal weight enumerators in $R_{4 \pm 2\sqrt{2}}^-$. The following theorem shows that there exists an extremal formal weight enumerator which does not satisfy the Riemann hypothesis: 
%%%%%%%%%%%%%%%%%% theorem %%%%%%%%%%%%%%%%%%
\begin{thm}\label{thm:4+2sqrt2_deg8}
Let $P_8^+(T)$ and $P_8^-(T)$ be the zeta polynomials of $\varphi_8^+ (x,y)$ and $\varphi_8^- (x,y)$, respectively. Then we have the following:\\
(i) Some of the roots $\alpha$ of $P_8^+(T)$ do not satisfy $|\alpha|=1/\sqrt{4 + 2\sqrt{2}}$ (the Riemann hypothesis is false). \\
(ii) All the roots $\alpha$ of $P_8^-(T)$ satisfy $|\alpha|=1/\sqrt{4 - 2\sqrt{2}}$ (the Riemann hypothesis is true). 
\end{thm}
\prf (i) Let $q=4 + 2\sqrt{2}$. We normalize the circle $|T|=1/\sqrt{q}$ to the unit circle $|T|=1$, so we consider  $P_8^+(T/\sqrt{q})$ instead of $P_8^+(T)$. We can calculate (with some computer algebra system) that 
$$P_8^+\left(\frac{T}{\sqrt{q}}\right)=(T-1)(T+1)(T^4-{\textstyle\sqrt{2+\sqrt{2}}} (T^3+T)-\sqrt{2}T^2+1).$$
The factor of degree four becomes 
$$T^2\left\{\left(T^2+\frac{1}{T^2}\right)-{\textstyle\sqrt{2+\sqrt{2}}}\left( T+\frac{1}{T} \right)-\sqrt{2}\right\}.$$
We put $V=T+1/T$. It suffices to solve 
$$V^2-{\textstyle\sqrt{2+\sqrt{2}}}V-(2+\sqrt{2})=0$$
and we get
$$V=\frac{(1\pm \sqrt{5})\sqrt{2+\sqrt{2}}}{2}.$$
It is easy to see that 
$$T^2-\frac{(1 + \sqrt{5})\sqrt{2+\sqrt{2}}}{2}T+1$$
has two distinct real roots which are not $\pm 1$. Thus we can conclude that two of the roots of $P_8^+(T)$ are not located on the circle $|T|=1/\sqrt{q}$ and that the Riemann hypothesis is false in this case (other four roots of $P_8^+(T)$ are on $|T|=1/\sqrt{q}$). 

\medskip
\noindent(ii) Let $q=4 - 2\sqrt{2}$. We have 
$$P_8^-\left(\frac{T}{\sqrt{q}}\right)=T^2(T-1)(T+1)
\left\{\left(T+\frac{1}{T}\right)^2+{\textstyle\sqrt{2-\sqrt{2}}} \left(T+\frac{1}{T}\right)-(2-\sqrt{2})\right\}.$$
We can prove similarly to (i) that all the roots of $P_8^-(T)$ lie on the circle $|T|=1/\sqrt{q}$, so the Riemann hypothesis is true. \qed
%%%%%%%%%%%%%%%%%% example %%%%%%%%%%%%%%%%%%
\begin{exam}\label{exam:4+2sqrt2_deg24}\rm 
We consider the ring $R^-_{4+2\sqrt{2}}$ (see (\ref{eq:ring_fwe_4+2sqrt2})). In the case of degree 24, the possible pairs of $(l,m)$ in (\ref{eq:form_fwe_4+2sqrt2}) are $(l,m)=(8,0)$ and $(0,1)$. We have

\medskip
$\displaystyle W(x,y):=\frac{21}{16}W_{2, 4 + 2\sqrt{2}}(x,y)^8\varphi_8^+(x,y)-\frac{5}{16}\varphi_8^+(x,y)^3$
\begin{eqnarray}
&=& x^{24}-(16422+11592 \sqrt{2})x^{20}y^4 +(1020096+721280\sqrt{2})x^{18}y^6\nonumber\\
& &-(33004977+23338008\sqrt{2})x^{16}y^8+(519785280+367543680\sqrt{2})x^{14}y^{10}\nonumber\\
& &-(4102489300+2900898000\sqrt{2})x^{12}y^{12}+(17657398080+12485665920\sqrt{2})x^{10}y^{14}\nonumber\\
& &-(38087686257+26932061232\sqrt{2})x^8y^{16}+(39988783296+28276339840\sqrt{2})x^6y^{18}\nonumber\\
& &-(21850472742+15450617448\sqrt{2})x^4y^{20}+(768398401+543339720\sqrt{2})y^{24}. \label{eq:4+2sqrt2_deg24}
\end{eqnarray}
There is no other formal weight enumerator of degree 24 with $d\geq 4$ in $R_{4 + 2\sqrt{2}}^-$, so (\ref{eq:4+2sqrt2_deg24}) is extremal. Using the intermediate value theorem on the real axis, we can verify that the zeta polynomial of $W(x,y)$ has a real root other than $\pm 1/\sqrt{q}$ ($q=4+2\sqrt{2}$) and that (\ref{eq:4+2sqrt2_deg24}) does not satisfy the Riemann hypothesis. 
\end{exam}

\medskip
There seem to be similar structures in the cases where $q=2 \pm 2\sqrt{5}/5$ and $8 \pm 4\sqrt{3}$. The polynomials $\varphi_{10}^{\pm}(x,y)$ and $\varphi_{12}^{\pm}(x,y)$ themselves are extremal (see (\ref{eq:fwe_deg10}) and (\ref {eq:fwe_deg12})). We have the following theorems: 
%%%%%%%%%%%%%%%%%% theorem %%%%%%%%%%%%%%%%%%
\begin{thm}\label{thm:2pm2sqrt5over5}
(i) The extremal formal weight enumerator $\varphi_{10}^{+}(x,y)$ does not satisfy the Riemann hypothesis. \\
(ii) The extremal formal weight enumerator $\varphi_{10}^{-}(x,y)$ satisfies the Riemann hypothesis. 
\end{thm}
\prf We use the intermediate value theorem on the real axis to the zeta polynomial $P_{10}^+(T)$ of $\varphi_{10}^{+}(x,y)$ to prove that it has real roots other than $T=\pm 1/\sqrt{q}$ ($q=2 + 2\sqrt{5}/5$) for (i). For (ii), we put $T=\cos\theta$ (so $T+1/T=2\cos\theta$) in the normalized  zeta polynomial $P_{10}^-(T/\sqrt{q})$ of $\varphi_{10}^{-}(x,y)$ ($q=2 - 2\sqrt{5}/5$) and use the same theorem to prove that all the solutions of $P_{10}^-(T/\sqrt{q})=0$ are on the unit circle (we omit the detail). \qed
%%%%%%%%%%%%%%%%%% theorem %%%%%%%%%%%%%%%%%%
\begin{thm}\label{thm:8pm4sqrt3}
(i) The extremal formal weight enumerator $\varphi_{12}^{+}(x,y)$ does not satisfy the Riemann hypothesis. \\
(ii) The extremal formal weight enumerator $\varphi_{12}^{-}(x,y)$ satisfies the Riemann hypothesis. 
\end{thm}
\prf Similar to Theorem \ref{thm:2pm2sqrt5over5}. \qed
%%%%%%%%%%%%%%%%%%%%%%%%%%%%%%%%%%%%%%%%%%%%%%%%%%%%
%%%%%%%%%%%%%%%%%%%%%%%%%%%%%%%%%%%%%%%%%%%%%%%%%%%%
\section{Some remarks}\label{section:remarks}
We can propose several problems concerning the contents of the previous sections. 

\medskip
\noindent (I) {\it On the factorization and roots of the polynomial $|A(n,q)|$}

\medskip
Here are some examples of factorization of $|A(n,q)|$ (see the discussion that follows Corollary \ref{cor:moment_even_deg} for the definition of $A(n,q)$): 

\medskip
$|A(1,q)|=-4,$

$|A(2,q)|=8(q-2),$ 

$|A(3,q)|=16(q-2)(3q-4),$ 

$|A(4,q)|=32(q-2)(3q-4)(q^2-8q+8),$ 

$|A(5,q)|=-64(q-2)(3q-4)(q^2-8q+8)(5q^2-20q+16),$ 

$|A(6,q)|=128(q-2)^2(3q-4)(q^2-8q+8)(5q^2-20q+16)(q^2-16q+16),$ 

$|A(7,q)|=256(q-2)^2(3q-4)(q^2-8q+8)(5q^2-20q+16)(q^2-16q+16)(7q^3-56q^2+112q-64),$ 

$|A(8,q)|=2|A(7,q)|(q^4-32q^3+160q^2-256q+128),$ 

$|A(9,q)|=-2|A(8,q)|(3q-4)(3q^3-36q^2+96q-64),$ 

$|A(10,q)|=-2|A(9,q)|(q-2)(q^4-48q^3+304q^2-512q+256),$ 

$|A(11,q)|=2|A(10,q)|(11q^5-220q^4+1232q^3-2816q^2+2816q-1024),$ 

$|A(12,q)|=2|A(11,q)|(q^2-8q+8)(q^4-64q^3+320q^2-512q+256).$ 

\medskip
\noindent It is not easy to predict the factorization of $|A(n,q)|$ in general, but Professor Masakazu Yamagishi points out the following: 
%%%%%%%%%%%%%%%%%% conj %%%%%%%%%%%%%%%%%%
\begin{conj}\label{conj:chebyshev}
$$|A(n,q)|=2(-1)^n q^{n/2} T_n(q^{-1/2})|A(n-1,q)|\qquad (n\geq 2),$$
where $T_n(x)$ is the Chebyshev polynomial of the first kind. 
\end{conj}

\medskip
\noindent (II) {\it On the Riemann hypothesis for extremal formal weight enumerators}

\medskip
In Section \ref{section:extremal_no_rh}, we considered the rings $R_{4\pm 2\sqrt{2}}^-$, 
\begin{eqnarray}
R_{2 \pm 2\sqrt{5}/5}^- &=& {\bf C}[W_{2, 2 \pm 2\sqrt{5}/5}(x,y), \varphi_{10}^{\pm}(x,y)], \label{eq:ringfwe_4+2sqrt5over5}\\
R_{8 \pm 4\sqrt{3}}^- &=& {\bf C}[W_{2, 8 \pm 4\sqrt{3}}(x,y), \varphi_{12}^{\pm}(x,y)]. \label{eq:ringfwe_8+4sqrt3}
\end{eqnarray}
and found several phenomena which are different from the ones known before in connection to the coding theory. The author has not found formal weight enumerators which satisfy the Riemann hypothesis in $R_{4 + 2\sqrt{2}}^-$, $R_{2 + 2\sqrt{5}/5}^-$ and $R_{8 + 4\sqrt{3}}^-$. On the other hand, there seems to be plenty of formal weight enumerators satisfying it in $R_{4 - 2\sqrt{2}}^-$, $R_{2 - 2\sqrt{5}/5}^-$ and $R_{8 - 4\sqrt{3}}^-$. At present, there is no explanation for this difference. We would like to ask the following questions: 
%%%%%%%%%%%%%%%%%% problem %%%%%%%%%%%%%%%%%%
\begin{prob}\label{prob:rh_extremal}
(i) Can one prove analogs of the Mallows-Sloane bound in the above rings ? \\
(ii) Are there (extremal) formal weight enumerators which satisfy the Riemann hypothesis in the rings $R_{4 + 2\sqrt{2}}^-$, $R_{2 + 2\sqrt{5}/5}^-$ and $R_{8 + 4\sqrt{3}}^-$ ? \\
(iii) Explain the difference of the truth of the Riemann hypothesis between the formal weight enumerators for $q=4+2\sqrt{2}$ and $4-2\sqrt{2}$ (the cases $2\pm 2\sqrt{5}/5$ and $8\pm 4\sqrt{3}$ are similar). 
\end{prob}

\bigskip
\noindent{\it Acknowledgement. }This work was supported by JSPS KAKENHI Grant Number JP26400028. The author would like to express his sincere gratitude to  Professor Masakazu Yamagishi for pointing out Conjecture \ref{conj:chebyshev}. This work was established mainly during the author's stay at University of Strasbourg for the overseas research program of Kindai University. He would like to express his appreciation to  Professor Yann Bugeaud at University of Strasbourg for his hospitality and to Kindai University for giving him a chance of the program.

%\end{linenumbers}

\begin{thebibliography}{99}
\bibitem[1] {Ch1} K. Chinen, Zeta functions for formal weight enumerators and the extremal property, Proc. Japan Acad. 81 Ser. A. (2005), 168-173. 
\bibitem[2] {Ch2} K. Chinen, An abundance of invariant polynomials satisfying the Riemann hypothesis, Discrete Math. 308 (2008), 6426-6440. 
\bibitem[3] {Ch3} K. Chinen, Extremal invariant polynomials not satisfying the Riemann hypothesis, arXiv:1709.03389. 

\bibitem[4] {Ch4} K. Chinen, On some families of divisible formal weight enumerators and their zeta functions, arXiv:1709.03396. 
\bibitem[5] {CoSl} J. H. Conway, N. J. A. Sloane, Sphere Packings, Lattices and Groups, third ed., Springer Verlag, NewYork, 1999.
\bibitem[6] {Du1} I. Duursma, Weight distribution of geometric Goppa codes, Trans. Amer. Math. Soc. 351, No.9 (1999), 3609-3639.
\bibitem[7] {Du2} I. Duursma, From weight enumerators to zeta functions, Discrete Appl. Math. 111 (2001), 55-73.
\bibitem[8] {Du3} I. Duursma, A Riemann hypothesis analogue for self-dual codes, DIMACS series in Discrete Math. and Theoretical Computer Science 56 (2001), 115-124.
\bibitem[9] {Du4} I. Duursma, Extremal weight enumerators and ultraspherical polynomials, Discrete Math. 268, No.1-3 (2003), 103-127. 
\bibitem[10] {MaSl} F. J. MacWilliams, N. J. A. Sloane, The Theory of Error-Correcting Codes, North-Holland, Amsterdam, 1977.  
\bibitem[11] {Oz} M. Ozeki, On the notion of Jacobi polynomials for codes, Math. Proc. Camb. Phil. Soc. 121 (1997), 15-30. 
\bibitem[12] {Pl} V. S. Pless, Introduction to the Theory of Error-Correcting Codes, third ed., John Wiley \& Sons, NewYork, 1998. 
\bibitem[13] {Ra} E. M. Rains, Shadow bounds for self-dual codes, IEEE Trans. Inform. Theory 44(1) (1998), 134-139. 
\bibitem[14] {ShTo} G. C. Shephard, J. A. Todd, Finite unitary reflection groups, Canad. J. Math. 6 (1954), 274-304. 
\bibitem[15] {Sl} N. J. A. Sloane, Self-dual codes and lattices, in: Relations Between Combinatorics and Other Parts of Mathematics, Proceedings of Symposium on Pure Mathematics 34, Ohio State University, Columbus, OH, 1978, Amer. Math. Soc., Providence, RI, 1979, pp. 273-308. 
\end{thebibliography}
\end{document}